\magnification=\magstep1
\input harvmac 

\hfuzz 15pt
\parskip=7pt  \baselineskip=12pt

\font\male=cmr9
\font\sfont=cmbx10 scaled\magstep1    %large

%%% DEFINITIONS 

\def\newsubsec#1{\global\advance\subsecno by1\message{(\secsym\the\subsecno.
#1)} \ifnum\lastpenalty>9000\else\bigbreak\fi 
\noindent{\bf\secsym\the\subsecno. #1}\writetoca{\string\quad
{\secsym\the\subsecno.} {#1}}}

\def\ha{\textstyle{1\over 2}} 
\def\lg{\langle}
\def\rg{\rangle}

\def\ll{\left(}
\def\Ll{\Big[ } %{\left[\ }
\def\LL{\left\{}
\def\rr{\right)} 
\def\Rr{ \Big]}%{\ \right]} 
\def\RR{\right\}} 
\def\nt{\noindent} 
\def\nl{\hfil\break}
\def\np{\vfil\eject} 

\def\bbz{Z\!\!\!Z}
\def\bbc{C\kern-6.5pt I}

\def\bbn{I\!\!N}

\def\PR{{\it Proof: }}
\def\ugh{\cu_{g,h}}
\def\up{\cu'_{g,h}}
\def\ut{\tilde\ugh}
\def\upt{\tilde\up}
\def\gl{GL_{g,h}(2)}

\def\ca{{\cal A}}
\def\cc{{\cal C}}
\def\cd{{\cal D}}

\def\cu{{\cal U}}

\def\cz{{\cal Z}}
    
\def\d{\delta}  \def\ve{\varepsilon}   
\def\g{\gamma}

  \def\rra{\longrightarrow}
 \def\mt{\mapsto}

\def\tc{\tilde C}
\def\ta{\tilde a}
\def\td{\tilde d} 
\def\tg{{\tilde g}}
\def\th{{\tilde h}}

%%%%%%%%% REFERENCES 

\lref\ACC{B. Abdesselam, A. Chakrabarti and R. Chakrabarti, 
Mod. Phys. Lett. {\bf A11} (1996) 2883.} 

\lref\Abe{E. Abe, {\it Hopf Algebras}, 
Cambridge Tracts in Math., N 74,(Cambridge Univ. Press, 1980).} 

\lref\Ag{A. Aghamohammadi, Mod. Phys. Lett. {\bf A8} (1993)
2607.} 

\lref\Aiz{N. Aizawa, preprint OWUAM-020, q-alg/9701022.}

\lref\BH{A. Ballesteros and F.J. Herranz, 
J. Phys. {\bf A29} (1996) L311.} 

\lref\BHOS{A. Ballesteros, F.J. Herranz, M.A. del Olmo and M.
Santander, J. Phys. {\bf A28} (1995) 941-956.}

\lref\DMMZ{E.E. Demidov, Yu.I. Manin, E.E. Mukhin and
D.V. Zhdanovich, Progr. Theor. Phys. Suppl. {\bf
102} (1990) 203-218.} 

\lref\Dod{V.K. Dobrev, J. Math. Phys. {\bf 33} (1992) 3419.}  

\lref\Dobh{V.K. Dobrev, 
in: Proceedings of the 10th International Conference
`Problems of Quantum Field Theory', (Alushta, Crimea, Ukraine,
13-18.5.1996), eds. D. Shirkov, D. Kazakov and A. Vladimirov,
JINR E2-96-369, (Dubna, 1996) pp. 104-110.}

\lref\DPa{V.K. Dobrev and P. Parashar, J. Phys. {\bf A26} (1993) 
6991.} 

\lref\DPb{V.K. Dobrev and P. Parashar, Lett. Math. Phys. {\bf 29}
(1993) 259.} 

\lref\Dr{V.G. Drinfeld, Dokl. Akad. Nauk SSSR {\bf 283}
(1985) 1060-1064; ~ in: Proceedings ICM, (MSRI, Berkeley,
1986) pp. 798-820.}

\lref\FRT{ L.D. Faddeev, N. Yu. Reshetikhin and L.A. Takhtajan, Alg.
Anal. {\bf 1} (1989) 178-206 (in Russian) and in: {\it Algebraic Analysis},
Vol. No. 1 (Academic Press, 1988) pp. 129-139. }

\lref\Jeu{J. Van der Jeugt, q-alg/9703011.}  

\lref\Kar{V. Karimipour,  Lett. Math. Phys. {\bf 30} (1994) 87.} 

\lref\KSAA{M. Khorrami, A. Shariati, M.R. Abolhassani and 
A. Aghamohammadi, Mod. Phys. Lett. {\bf A10}  (1995)  873.} 

\lref\Kup{B.A. Kupershmidt, J. Phys. {\bf A25} (1992) L1239.}   

\lref\Ma{Yu.I. Manin, Ann. Inst. Fourier {\bf 37} (1987) 191-205;
\ Quantum groups and non-commutative geometry, Montreal
University preprint, CRM-1561 (1988).} 

\lref\Ohn{C. Ohn, Lett. Math. Phys. {\bf 25} (1992) 85-88.}

\lref\Par{P. Parashar, preprint SISSA 85/96/FM; q-alg/9606003.}

\lref\QL{Podles and S.L. Woronowicz, Comm. Math. Phys. {\bf 130}
(1990) 381-431;\nl U. Carow-Watamura, M. Schlieker, M. Scholl and
S.  Watamura, Zeit. f. Physik {\bf C48} (1990) 159-166.} 

\lref\SWZ{A. Schirrmacher, J. Wess and B. Zumino, Z. Phys. {\bf
C49} (1991) 317.}   

\lref\Su{A. Sudbery, in: Proceedings of the Workshop on Quantum Groups, 
Argonne National Lab (1990), eds. T.  Curtright, D. Fairlie and C.
Zachos, (World Sci, 1991) pp. 33-52.}

\lref\Vla{A.A. Vladimirov, Mod. Phys. Lett. {\bf A8} (1993)
2573.}  
 
\lref\Zak{S. Zakrzewski, Lett. Math. Phys. {\bf 22} (1991)
287-289.} 

%%%%%%%%%%%%%%%  START of paper 

\line{\hfill INRNE-TH/3/97 (1997)} 

\vskip 1truecm

\centerline{{\sfont Duality for the Jordanian Matrix Quantum Group
~GL$_{g,h}$(2)}}  

\vskip 1.5cm

\centerline{{\bf B.L. Aneva}, ~{\bf V.K. Dobrev}$^*$ 
~and ~{\bf S.G. Mihov}}\footnote{}{$^*$\ 
{current address:}~ International Center for Theoretical Physics, 
via Costiera 11, P.O. Box 586, 34100 Trieste, Italy;~
Fax: +39-40-224163, ~ 
e-mail: dobrev@ictp.trieste.it}     
\vskip 0.5cm 

\centerline{Bulgarian Academy of Sciences} 
\centerline{Institute of Nuclear Research and Nuclear Energy} 
\centerline{72 Tsarigradsko Chaussee, 1784 Sofia, Bulgaria}

\vskip 1cm

\centerline{{\bf Abstract}}

\midinsert\narrower\narrower{\male 
We find the Hopf algebra ~${\cal U}_{g,h}$~ dual to the Jordanian
matrix quantum group ~$GL_{g,h}(2)$. As an algebra it depends
only on the sum of the two parameters and is split in two
subalgebras: ~${\cal U}'_{g,h}$ (with three generators) and ~$U(\cz)$
(with one generator). The subalgebra $U(\cz)$ is a central Hopf
subalgebra of ${\cal U}_{g,h}$. The subalgebra ~${\cal U}'_{g,h}$~ is not a
Hopf subalgebra and its coalgebra structure depends on both
parameters. We discuss also two one-parameter special cases: $g
=h$ and $g=-h$. The subalgebra ${\cal U}'_{h,h}$ is a Hopf algebra and
coincides with the algebra introduced by Ohn as 
the dual of $SL_h(2)$. The subalgebra ${\cal U}'_{-h,h}$ is isomorphic
to $U(sl(2))$ as an algebra but has a nontrivial coalgebra
structure and again is not a Hopf subalgebra of ${\cal U}_{-h,h}$.  }
\endinsert 

\vskip 2cm 

\newsec{Introduction}

\nt
The group ~$GL(2)$~ admits two distinct quantum group deformations 
with central quantum determinant: ~$GL_q(2)$ \Dr\ and 
~$GL_h(2)$ \DMMZ, \Zak. These are the only possible such 
deformations (up to isomorphism) \Kup. Both may be viewed as
special cases of two parameter deformations: ~$GL_{p,q}(2)$ \DMMZ\ and
~$GL_{g,h}(2)$ \Ag. In the initial years of the development of quantum
group theory mostly $GL_q(2)$ and $GL_{p,q}(2)$ were considered. More
recently started research on ~$SL_h(2)$~ and its dual quantum
algebra ~$U_h(sl(2))$~ \Ohn. In particular, aspects of 
differential calculus \Ag, and differential geometry \Kar\ were
developed, the universal R-matrix for $U_h(sl(2))$ was given
in \Vla, \KSAA, \BH,  representations of $U_h(sl(2))$ were 
constructed in \Dobh, \ACC, \Aiz, \Jeu, contractions of 
$SL_h(2)$ and $U_h(sl(2))$ were given in \Par. However, there 
are no studies until now of the two-parameter Jordanian matrix 
quantum group ~$GL_{g,h}(2)$. 
Even the dual of this algebra is not known. 

This is the problem we solve in this paper. 
We find the Hopf algebra ~$\cu_{g,h}$~ dual to the Jordanian
matrix quantum group ~$GL_{g,h}(2)$. As an algebra it depends on 
the sum ~$\tg ~=~ (g+h)/2$~ of the two parameters and is split in two
subalgebras: ~$\cu'_{g,h}$ (with three generators) and ~$U(\cz)$ 
(with one generator). The subalgebra $U(\cz)$ is a central Hopf
subalgebra of $\cu_{g,h}$. The subalgebra ~$\cu'_{g,h}$~ is not a
Hopf subalgebra and its coalgebra structure depends on both 
parameters. We discuss also two interesting one-parameter special
cases: $g =h$ and $g=-h$. The subalgebra $\cu'_{h,h}$ is a Hopf
algebra and coincides with the algebra introduced by Ohn as the
dual of $SL_h(2)$. The subalgebra $\cu'_{-h,h}$ is isomorphic to
$U(sl(2))$ as an algebra but has a nontrivial coalgebra structure
and again is not a Hopf subalgebra of $\cu_{-h,h}$.  

The paper is organized as follows. In Section 2 we recall the
group $GL_{g,h}(2)$. In Section 3 we first recall the method
(developed by one of us) of finding the dual. We then make a
change of generators and introduce the appropriate PBW basis in
$GL_{g,h}(2)$. We find the dual algebra ~$\ugh$~ 
and we give explicitly its algebra and
coalgebra structure in Propositions 1 and 2. We make a further
change of basis in order to bring the algebra to a form closer to
Ohn's \Ohn. Our main result is summarized in a Theorem. We also
introduce a subalgebra ~$\ut$~ of $\ugh$ with a basis such that 
no exponents of generators appear explicitly in the algebra and 
coalgebra relations. In Section 4 we consider the two interesting
one-parameter special cases: $g =h$ and $g=-h$ recovering the
algebra of \Ohn\ in the first special case. 
In an Appendix we apply the nonlinear map of \ACC\ to our 
dual algebra. 

\bigskip

\newsec{Jordanian matrix quantum group ~GL$_{g,h}$(2)} 

\nt
In this Section we recall the Jordanian two parameter deformation
~$GL_{g,h}(2)$~ of ~$GL(2)$~ introduced in \Ag\ (and denoted
$GL_{h,h'}$). One starts with a
unital associative algebra generated by four elements ~$a,b,c,d$~
of a quantum matrix ~$M ~=~ \pmatrix{a & b \cr c &d\cr}$~ with
the following relations ($g,h \in \bbc$)~:
\eqn\glh{ 
\eqalign{ &[a,c] ~=~ gc^2 ~, \qquad 
[d,c] ~=~ hc^2 ~, \qquad [a,d] ~=~ gdc-hac 
\cr 
&[a,b] ~=~ h({\cal D} -a^2) ~, \qquad 
[d,b] ~=~ g({\cal D} - d^2) ~, \qquad 
[b,c] ~=~ gdc+hac - ghc^2 \cr 
&\cd ~=~ ad-bc+hac ~=~ ad-cb-gdc + ghc^2  \cr }} 
where ~$\cd$~ is a multiplicative quantum determinant which is
not central (unless $g=h$)~:
\eqn\dcd{ 
[ a,\cd ] ~=~ [ \cd , d ] ~=~ (g-h)\cd c ~, \qquad 
[b, \cd ] ~=~ (g-h)(\cd d - a \cd) ~, \qquad 
[ c, \cd] ~=~ 0 } 
Relations \glh\ are obtained by applying either 
the method of Faddeev, Reshetikhin and Takhtajan \FRT, namely, 
by solving the monodromy equation:
$$R~M_1~M_2 ~~=~~ M_2~M_1~R $$
($M_1 = M {\hat \otimes} I_2$, $M_2 = I_2 {\hat \otimes} M$),
with $R$-matrix:  
\eqn\rrr{ 
R ~~=~~ \pmatrix{ 
1 & -h & h&gh \cr 
0 & \ 1 & 0 &-g \cr
0 & \ 0 & 1 &\ g \cr
0 &\  0 & 0 &\ 1 \cr}}
or the method of Manin \Ma\ using $M$ as transformation matrix of
the appropriate quantum planes \Ag.

The above algebra is turned into a  bialgebra ~$A_{g,h}(2)$~ with
the standard $GL(2)$ co-product ~$\d$~ and co-unit ~$\ve$~: 
\eqn\copr{\d \left( \pmatrix{a & b \cr c &d\cr} \right) ~=~ 
\pmatrix{a \otimes a + b
\otimes c & a \otimes b + b \otimes d \cr 
c \otimes a + d \otimes c & c \otimes b + d \otimes d\cr} }
\eqn\cou{ 
\ve \left( \pmatrix{a & b \cr c &d\cr} \right) ~=~ \pmatrix{1 & 0
\cr 0 & 1 \cr} } 
From \copr, resp., \cou\ it follows: 
\eqn\ddd{
\d (\cd) ~=~ \cd\otimes \cd ~, \qquad \ve (\cd) ~=~ 1 } 

Further, we shall suppose that $\cd$ is invertible, i.e., there
is an element ~$\cd^{-1}$~ which obeys:
\eqn\dud{ \cd \cd^{-1} ~=~ \cd^{-1} \cd ~=~ 1_\ca ~, \qquad 
(\cd^{-1}) ~=~ \cd^{-1}\otimes \cd^{-1} ~, \qquad \ve
(\cd^{-1}) ~=~ 1 }  
(Alternatively one may say that the algebra is extended with the
element ~$\cd^{-1}$.)  
In this case one defines the left and right inverse matrix of $M$
\Ag: 
\eqn\inv{ M^{-1} ~=~ \cd^{-1} \pmatrix{d+gc& -b+g(d-a)+g^{2}c 
\cr -c &a-gc\cr} ~=~ \pmatrix{d+hc & -b+h(d-a)+h^{2}c 
\cr -c&  a-hc \cr} \cd^{-1} }

The quantum group ~$GL_{g,h}(2)$~ is ~{\it defined}~ as the 
Hopf algebra obtained from the bialgebra ~$A_{g,h}(2)$~ when 
~$\cd^{-1}$~ exists and with antipode given by the formula: 
\eqn\anto{ \g(M) ~=~ M^{-1}  \qquad \Rightarrow \qquad 
\g(\cd) ~=~ \cd^{-1} ~, \quad \g(\cd^{-1}) ~=~ \cd }

For ~$g=h$~ one obtains from ~$GL_{g,h}(2)$~ the matrix quantum
group ~$GL_h(2) ~=~ GL_{h,h}(2)$, and, if the condition  ~$\cd
~=~ 1_\ca$~ holds, the matrix quantum group ~$SL_h(2)$.
Analogously, for ~$g=h=0$~ one obtains from ~$GL_{g,h}(2)$~ the
algebra of functions over the classical groups ~$GL(2)$~ and
~$SL(2)$, resp. 

\newsec{The dual of GL$_{g,h}$(2)} 

\newsubsec{Summary of the method} 

\nt 
Two bialgebras ~$\cu , \ca$~ are said to be ~{\it in duality}~
\Abe\ if there exists a doubly nondegenerate bilinear form 
\eqn\dua{\lg ~, ~\rg ~:~ \cu \times \ca \rra \bbc ~, ~~~
\lg ~, ~ \rg ~:~ (u,a) ~\mt ~\lg ~u ~, ~a ~ \rg
~, ~u \in\cu ~, ~a \in \ca }
such that, for ~$u,v\in\cu~, a,b\in\ca$~:
\eqna\dub
$$\eqalignno{&\lg ~ u ~, ~ab ~ \rg ~=~ \lg ~ \d_\cu(u) ~,
~a \otimes b ~ \rg ~, ~~~\lg ~ uv ~, ~a ~ \rg ~=~
\lg ~ u \otimes v ~,  ~\d_\ca (a) ~ \rg 
&\dub a\cr 
&\lg 1_\cu , a \rg ~=~ \ve_\ca (a) ~, ~~~\lg u ,
1_\ca \rg ~=~ \ve_\cu (u) &\dub b\cr }$$ 
Two Hopf algebras ~$\cu , \ca$~ are said to be ~{\it in duality}~
\Abe\ if they are in duality as bialgebras and if 
$$\eqalignno{&\lg \g_\cu (u) , a \rg ~=~ \lg u , \g_\ca (a)
\rg &\dub c\cr }$$ 

It is enough to define the pairing \dua\ between the generating
elements of the two algebras. The pairing between any other
elements of ~$\cu, ~\ca$~ follows then from relations \dub{} and
the standard bilinear form inherited by the tensor product. 

The duality between two bialgebras or Hopf algebras may be used
also to obtain the unknown dual of a known algebra. For that it
is enough to give the pairing between the generating elements of
the unknown algebra with arbitrary elements of the PBW basis of
the known algebra. Using these initial pairings and the duality
properties one may find the unknown algebra. 
Such an approach was first given by Sudbery \Su. He
obtained ~$U_q(sl(2)) \otimes U(u(1))$~ as the algebra of
tangent vectors at the identity of ~$GL_q(2)$. The initial pairings
were defined through the tangent vectors at
the identity. However, such calculations become very difficult
for more complicated algebras. Thus, in \Dod\ a
generalization was proposed in which the initial pairings are postulated to
be equal to the classical undeformed results. This generalized
method was applied in \Dod\ to the standard two-parameter
deformation $GL_{p,q}(2)$, (where also Sudbery's method was
used), then in \DPa\ to the multiparameter deformation of
$GL(n)$, and in \DPb\ to the matrix quantum Lorentz group of \QL. 
One should note that the dual of $GL_{p,q}(2)$ was obtained also
in \SWZ\ by methods of q-differential calculus. 

\newsubsec{Change of basis for GL$_{g,h}$(2) and generators of
the dual algebra} 

\nt
In the present paper we apply the method of \Dod\ to find the
dual of ~$GL_{g,h}(2)$. Following \Dod\ we first need to fix a PBW
basis of $\gl$. At first one may be inclined to use a PBW basis
as the one introduced in \Ohn\ for the case $g=h$, namely
consisting of all monomials ~$a^k d^\ell b^m c^n$, where ~$k,
\ell , m, n \in \bbz_+\,$. (Actually, the basis in \Ohn\ is for
$SL_h(2)$ and is obtained by restricting to indices fulfilling
$k\ell =0$.) However, the calculation with such a basis are more
difficult. Our analysis showed that it would be
simpler to work with the following PBW basis:
\eqn\bsa{a^k\ d^\ell\ c^n\ b^m ~, \qquad k, \ell, m, n \in \bbz_+ } 
the explanation being that the elements ~$a,d,c$~ generate a
subalgebra (though not a Hopf subalgebra) of $\gl$, cf. the first
line of \glh. Further simplification results if we make the
following change of generating elements and parameters: 
\eqn\cbs{\eqalign{
&\ta ~=~ \half (a+d) ~, 
\qquad \td ~=~ \half (a-d) \cr 
&\tg ~=~ \half (g+h) ~, 
\qquad \th ~=~ \half (g-h) \cr }} 

With these generating elements and parameters the algebra
relations become: 
\eqn\comu{\eqalign{
&c\ta ~=~ \ta c - \tg c^2 ~, \qquad c\td ~=~ \td c - \th c^2 ~,
\qquad  \td\ta ~=~ \ta\td - \tg \td c + \th \ta c \cr 
& b\ta ~=~ \ta   b  + \tg  c   b  - 2 \th  \ta   \td  
+ 2 \tg  \td^2  + (\tg^2  - \th^2)  \ta   c  
+  \tg ( \th^2  - \tg^2 ) c^2  \cr 
& b\td ~=~ \td   b - \th  c   b  + 2 \tg  \ta   \td  
- 2 \th \td^2   + (\th^2   - \tg^2)  \td   c  
+ \th (\tg^2   - \th^2 ) c^2  \cr 
& bc ~=~ cb + 2 \tg \ta c - 2\th\td c + (\th^2 - \tg^2) c^2 \cr 
&\cd ~=~ \ta^2 - \td^2 - cb + (\tg^2 - \th^2) c^2 
- \tg \ta c + \th\td c}} 
Note that these relations are written in anticipation of the
PBW basis: 
\eqn\bsb{ f ~=~ f_{k,\ell,m,n} ~=~ \ta^k\ \td^\ell\ c^n\ b^m ~, \qquad
k, \ell, m, n \in \bbz_+ }  
Note also that the relations in the subalgebras generated by
~$a,d,c$~ and ~$\ta,\td,c$~ are isomorphic under the change:
~$a\mt\ta$, ~$d\mt\td$, ~$c\mt c$, ~$g\mt\tg$, ~$h\mt\th$, cf.
the first lines in \glh\ and \comu. 

The coalgebra relations become: 
\eqna\capr \eqna\cau \eqna\aan 
$$\eqalignno{ 
&\d \left( \pmatrix{\ta & b \cr c &\td\cr} \right) ~=~
 \pmatrix{
\ta \otimes \ta + \td \otimes \td + \ha b\otimes c + \ha c\otimes b 
& \ta \otimes b + \td \otimes b + b \otimes \ta - b \otimes \td \cr &\cr
c \otimes \ta + c \otimes \td + \ta \otimes c - \td \otimes c 
& \ta \otimes \td + \td \otimes \ta + \ha b\otimes c - \ha c\otimes b 
\cr} \qquad  \qquad  &\capr {}\cr &&\cr 
&\ve \left( \pmatrix{\ta & b \cr c &\td\cr} \right) ~=~ 
\pmatrix{1 & 0 \cr 0 & 0 \cr}  &\cau {}\cr &&\cr 
&\eqalign{ 
\g \left( \pmatrix{\ta & b \cr c &\td\cr} \right) ~=&~ 
\cd^{-1} \pmatrix{   \ta -\td + (\tg+\th)c & 
-b - 2(\tg+\th)\td  +  (\tg+\th)^{2}c 
\cr -c &   \ta + \td - (\tg+\th)c \cr} ~= \cr 
=&~ \pmatrix{
  \ta -\td+(\tg-\th)c & 
 -b + 2(\th-\tg)\td+  (\tg-\th)^{2}c  
\cr -c&  \ta +\td + (\th-\tg)c 
\cr} \cd^{-1}}  &\aan {}\cr }$$ 

Let us denote by ~$\ugh ~=~ U_{g,h}(gl(2))$~ the unknown yet dual
algebra of $\gl$, 
and by ~$A,B,C,D$~ the four generators of ~$\ugh$. Following
\Dod\ we shall define the pairing ~$\lg Z, f\rg$, ~$Z=A,B,C,D$,
~$f$~ is from \bsb, as the classical tangent vector at the
identity: 
\eqn\duda{
\lg ~Z ~, ~f ~\rg ~~\equiv ~~ 
\ve \left( 
{\partial f \over \partial y} \right) ~, \qquad 
(Z,y) ~=~ (A,\ta),\ (B,b),\ (C,c),\ (D,\td) } 
From this we get the explicit expressions:
\eqna\dudb 
$$\eqalignno{
\lg ~A ~, ~f ~\rg ~~=~~ 
\ve\left({\partial f \over \partial 
\ta}\right) ~~=&~~  k  \d_{\ell 0} \d_{m0} \d_{n0}   &\dudb a\cr 
\lg ~B ~, ~f ~\rg ~~=~~ 
\ve\left({\partial f \over \partial 
b}\right) ~~=&~~ \d_{\ell 0} \d_{m1} \d_{n0}   &\dudb b\cr 
\lg ~C ~, ~f ~\rg ~~=~~ 
\ve\left({\partial f \over \partial 
c}\right) ~~=&~~ \d_{\ell 0} \d_{m0} \d_{n1}   &\dudb c\cr 
\lg ~D ~, ~f ~\rg ~~=~~ 
\ve\left({\partial f \over \partial 
\td}\right) ~~=&~~ \d_{\ell 1} \d_{m0} \d_{n0}   &\dudb d\cr }$$

\newsubsec{Algebra structure of the dual } 

\nt
First we find the commutation relations between the generators of
~$\ugh$. Below we shall need expressions like $e^{\nu B}$ which
we define as formal power series ~$e^{\nu B} ~=~ 1_\cu + 
\sum_{p\in\bbn} {\nu^p \over p!} B^p$. We have:  

{\bf Proposition 1:}~~ The commutation relations of the generators
~$A,B,C,D$~ introduced by \dudb{} are:
\eqna\cma  
$$\eqalignno{
 [B,C]  ~ =&~   D 
&\cma a\cr &&\cr
  [D,B]  ~ =&~  {\textstyle{1\over \tg}} (e^{2\tg B} - 1_\cu) 
    &\cma b\cr &&\cr
  [D,C]  ~ =&~  -2 C + \tg D^2 - \tg A
 &\cma c\cr &&\cr 
  [A,B]  ~ =&~  0    ~, \qquad 
  [A,C]  ~ =~  0    ~, \qquad 
  [A,D]  ~ =~  0 &\cma d\cr}$$
\PR  Using the assumed duality the above relations are
shown by calculating their pairings with the 
basis monomials ~$f ~=~ \ta^k \td^\ell c^n b^m$~ of the dual
algebra. In particular, the pairing of ~$f ~=~ 
\ta^k \td^\ell c^n b^m$~ with the commutators is: 
\eqna\cmmu 
$$\eqalignno{
\lg \ [B,C] \ , \ f \  \rg ~ =&~   \d_{\ell 1} \d_{m0} \d_{n0} 
&\cmmu a\cr 
\lg \  [D,B] \ , \ f \  \rg ~ =&~  \d_{\ell 0} \theta_{m1} \d_{n0}
2^m \tg^{m-1}   &\cmmu b\cr
\lg \  [D,C] \ , \ f \  \rg ~ =&~  - 2 \d_{\ell 0} \d_{m0} \d_{n1} +
2 \tg   \d_{\ell 2} \d_{m0} \d_{n0} 
&\cmmu c\cr 
\lg \  [A,B] \ , \ f \  \rg ~ =&~ \lg \  [A,C] \ , \ f \  \rg ~=~ 
\lg \  [A,D] \ , \ f \  \rg ~=~ 0 
&\cmmu d\cr 
&\theta_{rs} ~\equiv ~\cases{ 1 ~~&~~ $r\geq s$\cr 
0 ~~&~~ $r< s$\cr } &\cmmu e\cr
}$$
To calculate a commutator ~$\lg   [W,Z]  ,  f   \rg$~ one first 
calculates ~$\lg   WZ  ,  f   \rg$~ and ~$\lg  ZW  ,  f   \rg$. 
The pairing of any quadratic monomial of the
unknown dual algebra with ~$f ~=~ \ta^k \td^\ell c^n b^m$~ 
is given by the duality properties \dub{}:
\eqn\duc{\lg \ WZ \ , \ f \  \rg ~=~ \lg \ W\otimes Z \ , \ \d_\ca(f) \
\rg ~=~ \lg \ W\otimes Z \ , \ \sum_j f'_j \otimes f''_j\ \rg ~=~ 
\sum_j \lg \ W\ , \ f'_j \ \rg\  \lg \ Z\ , \ f''_j \ \rg } 
where ~$f'_j\,$, ~$f''_j$ are elements of the basis \bsb\ and so
further a direct application of \dudb{} is used. We should note that
these calculations though complicated do not require explicit
knowledge of ~$\d_\ca(f) = \sum_j f'_j \otimes f''_j $~ for all 
$f$, and furthermore not all terms in the sums are necessary. In 
particular, while calculating $\d_\ca(f)$ one may neglect terms
containing the element ~$c$~ on either side of the tensor sign in
second and higher degrees even before reordering the terms to the
basis monomials, since from the commutation relations it is clear
that those terms will not produce any term with $c$ in zero or 
first degree, and any $f'_j$ ($f''_j$) containing $c$ in
second and higher degrees will give zero in \duc\ for any $W$
($Z$). For the same reasons, if $W\neq C$ ($Z\neq C$) one may neglect 
terms containing the element ~$c$~ on the left (right) side of
the tensor sign in first degree even before reordering the terms
to the basis monomials. Similar reasons hold for the elements
$\td$. Taking into account such simplifications one may find
the pairings of the quadratic monomials necessary for \cmmu{},
e.g., 
\eqna\cmmv 
$$\eqalignno{
\lg \ BC \ , \ f \  \rg ~ =&~  \ha \d_{\ell 1} \d_{m0} \d_{n0} + \th 
\d_{\ell 1} \d_{m1} \d_{n0} + \d_{\ell 0} \d_{n0} \theta_{m2} 
\half(\tg^2-\th^2) \tg^{m-2} + \d_{\ell 0} \d_{m1} \d_{n1} \qquad
&\cmmv a\cr 
\lg \ CB \ , \ f \ \rg ~ =&~  -\ha \d_{\ell 1} \d_{m0} \d_{n0} + \th 
\d_{\ell 1} \d_{m1} \d_{n0} + \d_{\ell 0} \d_{n0} \theta_{m2} 
\half(\tg^2-\th^2) \tg^{m-2} + \d_{\ell 0} \d_{m1} \d_{n1} 
\qquad \qquad 
&\cmmv b\cr 
\lg \  DB \ , \ f \  \rg ~ =&~  \d_{\ell 0} \d_{n0} ( \d_{m1} + 
\theta_{m2} 2^{m-1} \tg^{m-2} (\tg-\th) ) 
&\cmmv c\cr
\lg \  BD \ , \ f \  \rg ~ =&~  - \d_{\ell 0} \d_{n0} ( \d_{m1} + 
\theta_{m2} 2^{m-1} \tg^{m-2} (\tg+\th) )
&\cmmv d\cr
\lg \  DC \ , \ f \  \rg ~ =&~  - \d_{\ell 0} \d_{m0} \d_{n1} +
( \th + \tg)   \d_{\ell 2} \d_{m0} \d_{n0} 
+ k \tg   \d_{\ell 1} \d_{m0} \d_{n0} 
+ \d_{\ell 1} \d_{m0} \d_{n1} 
&\cmmv e\cr
\lg \  CD \ , \ f \  \rg ~ =&~  \d_{\ell 0} \d_{m0} \d_{n1} +
( \th - \tg)   \d_{\ell 2} \d_{m0} \d_{n0} 
+ k \tg   \d_{\ell 1} \d_{m0} \d_{n0} 
+ \d_{\ell 1} \d_{m0} \d_{n1} 
 &\cmmv f\cr}$$
Note that quadratic relations \cmmv{} depend on both parameters,
while  the commutation relations \cmmu{}, which follow from
\cmmv{}, depend only on the parameter $\tg$.\nl 
Now in order to establish \cma{a} it is enough to compare the 
the RHS of \cmmu{a} and \dudb{d}.
Further, for relation \cma{b} we use \cmmu{b} and:
$$\eqalignno{
\lg \  B^p \ , \ f \  \rg ~ =&~ p!  \d_{\ell 0} \d_{mp} \d_{n0} 
&\cmmv h\cr }$$
(proved by induction) and its consequence:
$$\eqalignno{
\lg \   (e^{2\tg B} - 1_\cu)  \ , \ f \  \rg ~=&~ 
\sum_{p\in\bbn} {(2\tg)^p\over p!} 
\lg \ B^p  \ , \ f \  \rg ~=~ \sum_{p\in\bbn} {(2\tg)^p\over p!}  p!
\d_{\ell 0} \d_{mp} \d_{n0} 
~= &\cr =&~ (2 \tg)^m \d_{\ell 0} \theta_{m1} \d_{n0} 
&\cmmv i\cr }$$ 
To establish \cma{c} we compare the RHS of \cmmu{c} with the 
appropriate linear combination of the right-hand-sides of three
equations, namely \dudb{a}, \dudb{c} and  
$$\eqalignno{
\lg \  D^2 \ , \ f \  \rg ~ =&~ 2  \d_{\ell 2} \d_{m0} \d_{n0} 
+ k \d_{\ell 0} \d_{m0} \d_{n0} 
&\cmmv g\cr}$$ 
This finishes the Proof.~~$\bullet$  

Note that the commutation relations \cma{} depend only on the
parameter ~$\tg$~ and that the generator ~$A$~ is central.  
This is similar to the situation of the dual algebra $\cu_{p,q}$ of the
standard matrix quantum group $GL_{p,q}$ the commutation 
relations of which depend only on the combination $q' = \sqrt{pq}$ and
also one generator is central \SWZ, \Dod. Here the
central generator appears as a central extension but this is
fictitious since this may be corrected by a change of basis,
namely, by replacing the generator $C$ by a generator $\tc$~:
\eqn\nbc{ C ~~=~~ \tc ~-~ {\textstyle {\tg \over 2}} A } 
With this only \cma{c}  changes to:
$$\eqalignno{   [D,\tc]  ~ =&~  -2 \tc + \tg D^2  &\cma {c'}\cr}$$ 
Besides this change we shall make a change of generating elements
of $\ugh$ in order to bring the commutation relations to a  form
closer to the algebra of \Ohn. 
Thus, we make the following substitutions:
\eqna\subs
$$\eqalignno{ 
&D ~~=~~ e^{\mu B}\ H\ e^{\nu B} &\subs a\cr 
&C ~~=~~ e^{\mu' B}\ Y\ e^{\nu' B} ~-~ {\textstyle {\tg \over 2}}
\sinh ( \tg B) e^{(\mu'+\nu') B} 
~-~ {\textstyle {\tg \over 2}} A &\subs b\cr }$$
Substituting \subs{} into \cma{a} we get the desired result 
~$[B,Y] ~=~ H$~ if we choose: ~$\mu'=\mu$, ~$\nu'=\nu$. 
Substituting \subs{} into \cma{b} we get the desired result 
~$[H,B] ~=~ {2\over \tg} \sinh(\tg\ B)$~ if we choose: ~$
\mu+\nu = \tg$. Thus with conditions:
$$\eqalignno{ 
&\mu+\nu ~=~ \tg  ~, \qquad 
\mu' ~=~ \mu, \qquad \nu'=\nu  &\subs c\cr }$$ 
we obtain the following commutation relations instead of \cma{}~:
\eqna\cmb 
$$\eqalignno{ 
[B,Y] ~=&~ H 
&\cmb a\cr &&\cr 
[H,B] ~=&~ {\textstyle{2\over \tg}} \sinh(\tg B) 
&\cmb b\cr &&\cr
[H,Y] ~=&~ -  Y \cosh(\tg B) ~-~  \cosh(\tg B) Y ~=&\cr 
=&~ - 2 Y \cosh(\tg B) ~-~ \tg H \sinh(\tg B) 
~+~ \tg \sinh(\tg B) \cosh(\tg B) 
&\cmb c\cr &&\cr 
  [A,B]  ~ =&~  0    ~, \qquad 
  [A,Y]  ~ =~  0    ~, \qquad 
  [A,H]  ~ =~  0 &\cmb d}$$

Note that relations \cmb{a,b,c} coincide with those of the
one-parameter algebra of \Ohn, (cf. Subsection 4.1), 
though the coalgebra structure is different as we shall see
below. We can use this coincidence to derive the Casimir operator
of ~$\ugh$~: 
\eqn\cas{ \eqalign{
&\hat\cc_2 ~~=~~ f_1(A)\ \cc_2 ~+~ f_2(A) \cr 
&\cc_2 ~~=~~ \half\ (H^2 + \cosh^2 (\tg B)) ~+~
{\textstyle{1\over \tg}}\ (Y\sinh(\tg B) + \sinh(\tg B) Y) \cr}}
where ~$f_1(A)$, ~$f_2(A)$ are arbitrary polynomials in the central
generator $A$. To derive \cas\ it is enough to check that 
~$[\cc_2, Z] ~=~ 0$~ for ~$Z = B,Y,H$. The latter follows also from
the fact \BHOS\ that $\cc_2$ is the Casimir of the one-parameter
algebra of \Ohn. 

Finally we also write a subalgebra ~$\ut$~ of $\ugh$ with the 
basis~:  ~$A$, ~$K=e^{\tg B} = K^+$,  ~$K^{-1}=e^{-\tg B}
= K^-$, ~$Y$, ~$H$, so that in terms of $A$, $K$, $K^{-1}$, $Y$, $H$ 
no exponents of generators appear in the algebra and coalgebra relations. 
Thus instead of \cmb{} we have: 
\eqna\cmd  
$$\eqalignno{
 [K^\pm,Y]  ~ =&~ \pm\ \tg HK^\pm ~\pm\ {\textstyle{\tg \over 2}} 
 (1_\cu-K^{\pm 2})  &\cmd a\cr &&\cr
[H,K^{\pm}]  ~ =&~  K^{\pm 2} - 1_\cu &\cmd b\cr &&\cr
[H,Y] ~=&~ -  Y (K + K^{-1}) ~+~ {\textstyle{\tg \over 2}} H
(K^{-1} - K)  ~+~ {\textstyle{\tg \over 4}} (K^2 - K^{-2})
&\cmd c\cr &&\cr 
  K K^{-1} ~ =&~  K^{-1} K ~=~ 1_\cu   &\cmd {c'}\cr &&\cr 
  [A,K]  ~ =&~ [A,K^{-1}] ~=~ 0    ~, \qquad 
  [A,Y]  ~ =~  0    ~, \qquad 
  [A,H]  ~ =~  0 ~, \qquad &\cmd d\cr
  }$$ 

\newsubsec{Coalgebra structure of the dual }

\nt 
We turn now to the coalgebra structure of ~$\ugh$. We have: 

\nt
{\bf Proposition 2:}~~ 
(i) ~The comultiplication in the algebra ~$\ugh$~ is given by: 
\eqna\coa 
$$\eqalignno{
\d_\cu (A) ~=&~ A \otimes 1_\cu ~+~ 1_\cu \otimes A 
&\coa a\cr &&\cr 
\d_\cu (B) ~=&~ B \otimes 1_\cu ~+~ 1_\cu \otimes B 
&\coa b\cr &&\cr 
\d_\cu (Y) ~=&~ Y \otimes e^{-\tg B} ~+~ e^{\tg B} \otimes Y ~-~ 
{\textstyle{\th^2 \over \tg}}\ \sinh (\tg B) \otimes
A^2 e^{-\tg B} ~+~ \th H \otimes A e^{-\tg B} \qquad 
&\coa c\cr &&\cr 
\d_\cu (H) ~=&~ H \otimes e^{-\tg B} ~+~ e^{\tg B} \otimes H ~-~ 
{\textstyle{2\th\over\tg}} \ \sinh (\tg B) \otimes A e^{-\tg B} 
&\coa d\cr }$$ 
(ii) ~The co-unit relations in ~$\ugh$~ are given
by: 
\eqn\coug{\ve_\cu (Z) ~=~ 0 ~, \qquad Z ~=~ A, B, Y, H } 
(iii) ~The antipode in the algebra ~$\ugh$~ is given by: 
\eqna\ant 
$$\eqalignno{
\g_\cu (A) ~=&~ - A
&\ant a\cr &&\cr 
\g_\cu (B) ~=&~ - B 
&\ant b\cr &&\cr 
\g_\cu (Y) ~=&~ - e^{-\tg B} Y e^{\tg B} ~+~ 
{\textstyle{\th^2\over \tg}}\ \sinh (\tg B) A^2 ~+~
\th  e^{-\tg B} H A e^{\tg B}  
&\ant c\cr &&\cr 
\g_\cu (H) ~=&~ - e^{-\tg B} H e^{\tg B} ~-~ 
{\textstyle{2\th\over\tg}} \sinh (\tg B) A 
&\ant d\cr }$$ 
\PR 
(i) ~We use the duality property \dub{a}, namely we  have 
$$\lg ~ Z ~, ~f_1\ f_2 ~ \rg ~=~ \lg ~ \d_\cu(Z) ~, ~f_1
\otimes f_2 ~ \rg $$
for every generator $Z$ of $\ugh$ and for every ~$f_1, f_2\in
\gl$. Then we calculate separately the  LHS 
and RHS and comparing the results prove \coa{}. The check of
\coa{a,b} is easy. Instead of \coa{c,d} we first find the coproduct
of the original generators $C,D$~: 
$$\eqalignno{
\d_\cu (C) ~=&~ C \otimes 1_\cu ~+~ e^{2\tg B} \otimes C ~+~ 
{\textstyle{1\over 2\tg}}\ ( e^{2\tg B} - 1_\cu ) \otimes
(\tg^2 A - \th^2 A^2) ~+~ \th D \otimes A \qquad \qquad 
&\coa {c'}\cr 
\d_\cu (D) ~=&~ D \otimes 1_\cu ~+~ e^{2\tg B} \otimes D ~+~ 
{\textstyle{\th\over\tg}} ( 1_\cu  - e^{2\tg B} ) \otimes A 
&\coa {d'}\cr }$$ 
and then \coa{c,d} follow. To check \coa{c'} the following
choices are crucial: 
~$(f_1,f_2) = (\ta^k \td^\ell c^n b^m,1_\cu)$, $(b^m,c)$, $(b^m,a^k)$,
$(d^\ell,a^k)$. To check \coa{d'} the crucial choices are:
~$(f_1,f_2) = (\ta^k \td^\ell c^n b^m,1_\cu)$, $(b^m,d)$,
$(b^m,a^k)$. 
\nl  
(ii) ~Formulae \coug\ follow from ~$\ve_\cu (Z) ~=~  \lg ~ Z , 1_\ca ~
\rg$, cf. \dub{b}, and using the defining relations
\dudb{}.\nl 
(iii) ~Formulae \ant{} follow from \dub{c} or by using the
following Hopf algebra axiom \Abe: 
\eqn\hop{m \circ ({\rm id} \otimes \g_\cu) \circ \d_\cu ~=~ i
\circ \ve_\cu } 
where ~$m$~ is the usual product in the
algebra: ~$m(Z\otimes W) ~=~ ZW~, ~Z,W\in\cu$~ and ~$i$~ is the
natural embedding of ~$F$~ into ~$\cu$~: ~$i(\nu) ~=~ \nu 1_\cu~,
~\nu\in F$. This is applied in our case with ~$\cu \mt \ugh$,
$F=\bbc$, to the elements $A,B,Y,H$ and using \coa{} and
\coug.~~$\bullet$  

{\bf Corollary 1:}~~ For later reference we mention also the
coproduct and antipode of the intermediate generator $\tc$ and
the antipode of the initial generator $D$:
$$\eqalignno{
\d_\cu (\tc) ~=&~ \tc \otimes 1_\cu ~+~ e^{2\tg B} \otimes \tc ~-~ 
{\textstyle{\th^2\over 2\tg}}\ ( e^{2\tg B} - 1_\cu ) \otimes A^2
~+~ \th D \otimes A \qquad \qquad
&\coa {c''}\cr &&\cr   
\g_\cu (\tc) ~=&~ - e^{-2\tg B} \tc ~+~ 
{\textstyle{\th^2\over 2\tg}}\ ( 1_\cu - e^{-2\tg B} ) A^2 ~+~
\th  e^{-2\tg B} D A  
&\ant {c'}\cr &&\cr 
\g_\cu (D) ~=&~ - e^{-2\tg B} D ~+~ 
{\textstyle{\th\over\tg}} ( e^{-2\tg B} - 1_\cu ) A 
&\ant {d'}\cr }$$

{\bf Corollary 2:}~~ The coalgebra structure in the 
subalgebra ~$\ut$~ is given as follows:\nl 
(i) ~comultiplication~:  
\eqna \coak 
$$\eqalignno{
\d_\cu (A) ~=&~ A \otimes 1_\cu ~+~ 1_\cu \otimes A 
&\coak a\cr &&\cr 
\d_\cu (K^\pm) ~=&~ K^\pm \otimes K^\pm 
&\coak b\cr &&\cr 
\d_\cu (Y) ~=&~ Y \otimes K^{-1} ~+~ K \otimes Y ~-~ 
{\textstyle{\th^2 \over 2\tg}}\ (K - K^{-1}) \otimes
A^2 K^{-1} ~+~ \th H \otimes A K^{-1} \qquad 
&\coak c\cr &&\cr 
\d_\cu (H) ~=&~ H \otimes K^{-1} ~+~ K \otimes H ~+~  
{\textstyle{\th\over \tg}} \ (K^{-1} - K) 
\otimes A K^{-1} 
&\coak d\cr }$$ 
(ii) ~co-unit~: 
\eqn\couk{\ve_\cu (Z) ~=~ 0 ~, \quad Z ~=~ A, Y, H , \qquad 
\ve_\cu (Z) ~=~ 1 ~, \quad Z ~=~ K, K^{-1} } 
(iii) ~antipode~: 
\eqna\antk 
$$\eqalignno{
\g_\cu (A) ~=&~ - A
&\antk a\cr &&\cr 
\g_\cu (K^\pm) ~=&~ K^{\mp}  
&\antk b\cr &&\cr 
\g_\cu (Y) ~=&~ - K^{-1} Y K ~+~ 
{\textstyle{\th^2\over 2\tg}}\ (K - K^{-1}) A^2 ~+~
\th  K^{-1} H A K 
&\antk c\cr &&\cr 
\g_\cu (H) ~=&~ - K^{-1} H K ~+~ 
{\textstyle{\th\over\tg}}\ (K^{-1} - K) A 
&\antk d\cr }$$ 

\newsubsec{Main result} 

\nt 
Finally we can state the following: 

\nt 
{\bf Theorem:}~~~ 
The Hopf algebra ~$\ugh$~ dual to ~$GL_{g,h}(2)$~ is generated by
~$A,B,Y,H$ (or  $A, K, K^{-1}, Y, H$), cf. relations \dudb{} and
\subs{}.  It is given by relations  \cmb{}, \coa{}, \coug, \ant{},
(resp.  \cmd{}, \coak{}, \couk, \antk{}). As an algebra it
depends only on one parameter ~$\tg ~=~ (g+h)/2$~ and is split in
two subalgebras:~ $\up$ (resp.  $\upt$) generated by ~$B,Y,H$
(resp. $K, K^{-1}, Y, H$) and ~$U(\cz)$, where the algebra
~$\cz$~ is spanned by ~$A$. The subalgebra $U(\cz)$ is central in
$\ugh$ and is also a Hopf subalgebra of ~$\ugh$. The subalgebra
~$\up$ (resp. $\upt$) is not a Hopf subalgebra.\nl
\PR 
Actually this statement is summarizing our results in this
section, cf. Propositions 1 and 2, and the basis change \subs{}. It
remains only to note that ~$U(\cz)$~ is a Hopf subalgebra since 
$A$ commutes with the other generators and its Hopf algebra
operations are in terms of $A$ itself. The subalgebra
generated by ~$\up$ (resp. $\upt$) 
is not a Hopf subalgebra since the 
generator ~$A$~ takes part in formulae \coa{c,d}, \ant{c,d} 
(resp. \coak{c,d}, \antk{c,d}).~~$\bullet$ 

\newsec{One-parameter cases} 

\nt 
It is interesting to discuss the one-parameter special cases of
the matrix quantum group $\gl$ and its dual.  

\newsubsec{Case g=h} 

\nt 
The one-parameter matrix quantum group ~$GL_\tg(2)$~ \DMMZ, \Zak, is
obtained from $\gl$ by setting $g=h=\tg$. Thus the dual algebra 
~$\cu_\tg ~\equiv ~\cu_{\tg,\tg}$~ of $GL_\tg(2)$ 
is obtained by setting ~$\th ~=~ \ha(g-h) ~=~ 0$~ in \cmb{}, \coa{},
\coug, \ant{}. Since the commutation relations \cmb{} and the
counit relations \coug\ do not depend on $\th$ they remain
unchanged for ~$\cu_\tg$. The coproduct and antipode relations of
~$\cu_\tg$~ are: 
\eqna\coag \eqna\antg 
$$\eqalignno{
\d_\cu (A) ~=&~ A \otimes 1_\cu ~+~ 1_\cu \otimes A 
&\coag a\cr &&\cr 
\d_\cu (B) ~=&~ B \otimes 1_\cu ~+~ 1_\cu \otimes B 
&\coag b\cr &&\cr 
\d_\cu (Y) ~=&~ Y \otimes e^{-\tg B} ~+~ e^{\tg B} \otimes Y 
&\coag c\cr &&\cr 
\d_\cu (H) ~=&~ H \otimes e^{-\tg B} ~+~ e^{\tg B} \otimes H 
&\coag d\cr &&\cr 
\g_\cu (A) ~=&~ - A
&\antg a\cr &&\cr 
\g_\cu (B) ~=&~ - B 
&\antg b\cr &&\cr 
\g_\cu (Y) ~=&~ - e^{-\tg B} Y e^{\tg B}
&\antg c\cr &&\cr 
\g_\cu (H) ~=&~ - e^{-\tg B} H e^{\tg B}
&\antg d\cr }$$ 

We see that the one-parameter Hopf algebra ~$\cu_\tg$~ is split
in two Hopf subalgebras ~$\cu'_\tg ~\equiv~ \cu'_{\tg,\tg}$~ and
$U(\cz)$ and we may write: 
\eqn\iss{  \cu_\tg ~~=~~ \cu'_\tg \otimes U(\cz) } 

Now we compare the algebra ~$\cu'_\tg$~ with the algebra of
\Ohn. We see that after the identification ~$B\mt X$, ~$\tg \mt
-h$, the algebra relations \cmb{a,b,c} and the 
coalgebra relations \coag{b,c,d}, \coug, \antg{b,c,d} coincide
with their counterparts in \Ohn, i.e., the algebra ~$\cu'_\tg$~
coincides with the algebra of Ohn. We also note that the algebra
~$\cu'_\tg$~ in the basis ~$B,\tc,D$~ (cf. \cma{a,b,c'},
\coa{b,c'',d'}, \coug, \ant{b,c',d'}) coincides for $\th=0$ with
the version given in \BH\ after the identification: 
~$(B,\tc,D;\tg) ~\mt ~ (A_+,A_-,A;z)$, and by using the opposite
coalgebra structure. 

\newsubsec{Case g = -h} 

\nt 
Here we consider another one-parameter case: ~$g = -h = \th$, i.e.,
~$\tg =0$. From \cmb{}, \coa{}, \ant{}, we obtain: 
\eqna\cmh  
$$\eqalignno{
[B,Y] ~=&~ H &\cmh a\cr &&\cr
[H,B]  ~=&~ 2 B     &\cmh b\cr &&\cr
 [H,Y] ~=&~ - 2 Y  &\cmh c\cr &&\cr 
  [A,B]  ~ =&~  0    ~, \qquad 
  [A,Y]  ~ =~  0    ~, \qquad 
  [A,H]  ~ =~  0 &\cmh d\cr}$$
\eqna\coah \eqna\anth 
$$\eqalignno{
\d_\cu (A) ~=&~ A \otimes 1_\cu ~+~ 1_\cu \otimes A 
&\coah a\cr &&\cr 
\d_\cu (B) ~=&~ B \otimes 1_\cu ~+~ 1_\cu \otimes B 
&\coah b\cr &&\cr 
\d_\cu (Y) ~=&~ Y \otimes 1_\cu ~+~ 1_\cu \otimes Y ~-~ 
\th^2 \ B \otimes A^2 ~+~ \th H \otimes A 
&\coah c\cr &&\cr 
\d_\cu (H) ~=&~ H \otimes 1_\cu ~+~ 1_\cu \otimes H ~-~ 2\th \ B
\otimes A  
&\coah d\cr &&\cr  
\g_\cu (A) ~=&~ - A
&\anth a\cr &&\cr 
\g_\cu (B) ~=&~ - B 
&\anth b\cr &&\cr 
\g_\cu (Y) ~=&~ - Y ~+~ \th^2 \ B A^2 ~+~ \th H A 
&\anth c\cr &&\cr 
\g_\cu (H) ~=&~ -  H ~-~ 2\th B A 
&\anth d\cr }$$ 

Thus, for ~$\tg ~=~ 0$~ the interesting feature 
is that the subalgebra ~$\cu'_{\th,-\th}$~ is isomorphic to the
undeformed ~$U(sl(2))$~ with $sl(2)$ spanned by $B,Y,H$. 
However, as in the general case, the 
coalgebra sector is not classical, and the generators $B,Y,H$ do
not close a co-subalgebra. 

\vskip 10mm

\nt 
{\bf Acknowledgments.}~~ B.L.A. was supported in part by BNFR under 
contract Ph-404, V.K.D. was supported in part by BNFR under 
contract Ph-643. 

\vskip 10mm

\parskip=0pt  \baselineskip=10pt
\appendix{A}{Application of a nonlinear map}%\foot{Added after submission.}} 

In \ACC\ a nonlinear map was proposed under which the one-parameter 
Ohn's algebra was brought to undeformed $sl(2)$ form, 
though, the coalgebra structure becomes even more complicated, 
cf. \Aiz\ and \Jeu. Since our two-parameter dual 
is like Ohn's algebra in the algebra sector we can also apply 
the map of \ACC, which we do in this Appendix.  
We give the map in our notation, namely, following (28) and 
(33) of \ACC\ we set:
\eqna\mpp
$$\eqalignno{
I_+ ~~=&~~ 
{2\over \tg}\ \tanh \left( {\tg\ B\over 2}\right) ~=~ 
- {2\over \tg}\ \left( 
1_\cu \ + \ 2\ \sum_{\ell=1}^{\infty}\ (-K)^\ell \right) 
~\left( ~=~ {2\over \tg}\ \left( 
{ K\ -\ 1_\cu\over K\ +\ 1_\cu }\right) \right) 
&\mpp a\cr  
&&\cr 
I_- ~~=&~~ \cosh \left( {\tg\ B\over 2}\right) 
\ Y \ \cosh \left( {\tg\ B\over 2}\right) 
~=~ \textstyle{1\over 4}\ \left(K^{1/2} + K^{-1/2}\right)\ Y\ 
\left(K^{1/2} + K^{-1/2}\right) \qquad &\mpp b\cr  
}$$
Then we have, following \ACC, (note though that   
we do not rescale $H$) 
the classical $sl(2)$ commutation relations and Casimir:
\eqn\cls{ [H,I_\pm] ~=~ \pm\ 2 \ I_\pm ~, \qquad [I_+,I_-] ~=~ H 
~, \qquad 
\cc_2^c ~=~ I_+ I_- \ + \ I_- I_+ \ +\ \half\ H^2 }

Of course, our aim is to write the coproducts. Actually, for $I^+$ 
we use (4.5) of \Jeu\ (since $I^+$ is expressed 
through $B$ which has the (parameter-independent) 
classical coproduct \coa{b} as in the one-parameter case) 
which in our notation gives: 
\eqn\coas{
\d_\cu (I_+) ~=~ I_+\ \otimes\ 1_\cu \ +\ 1_\cu\ \otimes\ I_+ \ +\ 
\sum_{n=1}^{\infty}\ 
\left( - { \tg^2 \over 4 }\right)^{n} \ 
\left( I_+^{n+1}\ \otimes\ I_+^{n}\ +\ 
I_+^{n}\ \otimes\ I_+^{n+1} \right) } 

For the co-product of $H$ we need the inverse of \mpp{a} 
(cf. (3.1) of \Aiz): 
\eqn\aiza{ 
 K^{\pm 1} ~=~ e^{\pm \tg B} ~=~ 1_\cu \ +\ 
2\ \sum_{\ell=1}^{\infty}\ \left( \pm {\tg \over 2}\ I_+ \right)^\ell 
~\left( ~=~ { 1_\cu \pm {\tg \over 2} I_+ 
\over 1_\cu \mp {\tg \over 2} I_+ } \right) }
Then we have using \coa{d}:
\eqn\coat{\eqalign{ 
\d_\cu (H) ~=&~ H \otimes 1_\cu\ +\ 1_\cu \otimes H \ +\ 
2\ \sum_{n=1}^{\infty} \left( H\ \otimes \ 
\left( -{\tg \over 2}\ I_+ \right)^n 
\ +\ \left( {\tg \over 2}\ I_+ \right)^n\ \otimes\ H \right)
~-\cr & \cr &-~ 
2\th  \ I_+\ 
\sum_{k=0}^{\infty}\ \left(  {\tg \over 2}\ I_+ \right)^{2k}
\otimes\ A\  
\left( 1_\cu \ +\ 
2\ \sum_{\ell=1}^{\infty}\ \left( - {\tg \over 2}\ I_+ \right)^\ell 
\right)}}

For the coproduct of $I_-$ we use \coa{c} and: 
\eqna\vco 
$$\eqalignno{
&\d_\cu (I_-) ~=~  
\d_\cu \left(\cosh \left( {\tg\ B\over 2}\right) \right)\  
\d_\cu (Y) \ 
\d_\cu \left(\cosh \left( {\tg\ B\over 2}\right) \right) &\vco a\cr 
&&\cr
&\d_\cu \left(\cosh \left( {\tg\ B\over 2}\right) \right) 
~=~ 
\cosh \left( {\tg\ B\over 2}\right) \otimes 
\cosh \left( {\tg\ B\over 2}\right) \ +\ 
\sinh \left( {\tg\ B\over 2}\right) \otimes 
\sinh \left( {\tg\ B\over 2}\right) \qquad \qquad 
&\vco b\cr 
}
$$
to obtain: 
\eqna\coay
$$ \eqalignno{
\d_\cu (I_-) ~=&~  \ I_- \ \otimes\  
\sum_{\ell=0}^{\infty}\ \ll \ell +1\rr\ 
\left( - {\tg \over 2}\ I_+ \right)^{\ell}  
 ~+~  \sum_{\ell=0}^{\infty}\ \ll \ell +1\rr\ 
\left(  {\tg \over 2}\ I_+ \right)^{\ell}  
\ \otimes\  I_- 
\ -  \cr & -\  
{\tg\over2}\ 
\ll I_+ \ I_-   \ + \ I_+ \ I_-   \rr \ 
\otimes\  
\sum_{\ell=1}^{\infty}\ \ell\ 
\left(  - {\tg \over 2}\ I_+ \right)^{\ell} 
\ +  \cr & +\  
{\tg\over2}\ 
\sum_{\ell=1}^{\infty}\ \ell\ 
\left(   {\tg \over 2}\ I_+ \right)^{\ell} 
\ \otimes\  
\ll I_+ \ I_-   \ + \ I_+ \ I_-   \rr 
\ +  \cr & +\  
{\tg^2\over4}\ I_+ \ I_-   \ I_+
\ \otimes\  
\sum_{\ell=2}^{\infty}\ \ll \ell -1\rr \ 
\left( - {\tg \over 2}\ I_+ \right)^{\ell} 
\ +  \cr & +\  
{\tg^2\over4}\ \sum_{\ell=2}^{\infty}\ \ll \ell -1\rr \ 
\left(  {\tg \over 2}\ I_+ \right)^{\ell} 
\ \otimes\  I_+ \ I_-   \ I_+
\ -\cr & \cr &- \  \th^2 \  \LL \  I_+ \ \otimes\ A^2\ \RR  
\ \left\{ \  \sum_{k=0}^{\infty}\ \ll k+1\rr 
\ \left(  {\tg \over 2}\ I_+ \right)^{2k}
\ \otimes\ 1_\cu
\ +  \right.  \cr &  \left.  + \ 
\sum_{k=0}^{\infty}\ \left(  {\tg \over 2}\ I_+ \right)^{2k} \ 
\otimes\  
 \sum_{\ell=1}^{\infty}\ 
\left( - {\tg \over 2}\ I_+ \right)^{\ell} 
\ +  \right.  \cr &  \left.  + \ 
\sum_{k=0}^{\infty}\ 
\ll k+1\rr \ \left( - {\tg \over 2}\ I_+ \right)^{k} 
\ \otimes\  
\sum_{\ell=1}^{\infty}\ \ell \ 
\left( - {\tg \over 2}\ I_+ \right)^{\ell} 
\ 
\right\}
\ +\cr & \cr &+ \ \th \  
\LL\ 1_\cu\ \otimes \ A\ \RR \ 
\Big\{ \  \Ll  H \ \otimes\ 1_\cu \Rr 
\ \times  \cr & \times   
\ \Ll  \ 
\sum_{k=0}^{\infty}\ \left(  {\tg \over 2}\ I_+ \right)^{2k} 
\ \otimes\  1_\cu 
\ +\ 
1_\cu \ \otimes\ 
\sum_{\ell=1}^{\infty}\ 
\ll \ell+1\rr \ \left( - {\tg \over 2}\ I_+ \right)^{\ell} 
\ +  \cr &   + \ 2\   
\sum_{k=1}^{\infty}\ \left( - {\tg \over 2}\ I_+ \right)^{k} 
\ \otimes\    
\sum_{\ell=1}^{\infty}\ \ell\ 
\left( - {\tg \over 2}\ I_+ \right)^{\ell} \ 
\Rr     
\ - &\coay {}\cr  &  -\     
2\ \Ll   \ 
\sum_{k=1}^{\infty}\ k\ 
\left( - {\tg \over 2}\ I_+ \right)^{2k} 
\ \otimes\ 1_\cu 
\ +\ \sum_{k=1}^{\infty}\ k\ 
\left( - {\tg \over 2}\ I_+ \right)^{k} 
\ \otimes\  
\sum_{\ell=1}^{\infty}\ \ell\ \left( 
- {\tg \over 2}\ I_+ \right)^{\ell}  
\ \Rr     
\Big\}
}$$ 

In the special case ~$\th=0$~ the coproducts of $H$ and $I_-$ 
coincide with the one-parameter formulae of \Aiz, cf. (3.2) and (5.3), 
resp., (with $\tg \mt -h$).   
In the special case ~$\tg=0$~ the nonlinear map becomes an identity 
and naturally the coproducts of $I_+$, $I_-$, $H$, concide 
with those of $B$, $Y$, $H$, resp., cf. \coah{b,c,d}. 

\vskip 10mm

\parskip=0pt 
\listrefs 

\np\end